\documentclass[11pt]{article}
\usepackage{amsmath}
\usepackage{amssymb,amsfonts,amsmath,latexsym,cite,color, psfrag}
\newtheorem{theo}{Theorem}

\newtheorem{prop}[theo]{Proposition}
\newtheorem{conj}[theo]{Conjecture}

\allowdisplaybreaks

\makeatletter \@addtoreset{equation}{section}
\@addtoreset{theo}{}\makeatother

\numberwithin{table}{section}

\def\qed{\hfill \rule{4pt}{7pt}}
\def\pf{\noindent {\it Proof.} }

\begin{document}
\begin{center}
{\Large\bf On  parabolic Kazhdan-Lusztig $R$-polynomials
for

 the symmetric group}

\vskip 6mm
{\small  Neil J.Y. Fan$^1$, Peter L. Guo$^2$, Grace  L.D. Zhang$^3$}

\vskip 4mm
$^1$Department of Mathematics\\
Sichuan University, Chengdu, Sichuan 610064, P.R. China
\\[3mm]

$^{2,3}$Center for Combinatorics, LPMC-TJKLC\\
Nankai University,
Tianjin 300071,
P.R. China

\vskip 4mm

$^1$fan@scu.edu.cn, $^2$lguo@nankai.edu.cn,
$^3$zhld@mail.nankai.edu.cn
\end{center}

\begin{abstract}
Parabolic  $R$-polynomials were introduced by Deodhar as
parabolic analogues of ordinary $R$-polynomials defined
 by Kazhdan and Lusztig. In this paper, we are
concerned with the computation of
parabolic  $R$-polynomials for the symmetric group.
Let $S_n$ be the symmetric group
on $\{1,2,\ldots,n\}$, and let
$S=\{s_i\,|\, 1\leq i\leq n-1\}$ be the generating
 set of $S_n$, where for $1\leq i\leq n-1$, $s_i$
 is the adjacent transposition. For a
  subset  $J\subseteq S$, let $(S_n)_J$ be the  parabolic subgroup
generated by $J$, and let
$(S_n)^{J}$
 be the set of minimal  coset representatives for $S_n/(S_n)_J$.
For $u\leq v\in (S_n)^J$ in the Bruhat order and $x\in \{q,-1\}$,
let $R_{u,v}^{J,x}(q)$ denote the parabolic  $R$-polynomial indexed by $u$ and $v$.
Brenti found a  formula for
$R_{u,v}^{J,x}(q)$ when   $J=S\setminus\{s_i\}$, and obtained an expression for
$R_{u,v}^{J,x}(q)$ when  $J=S\setminus\{s_{i-1},s_i\}$.
We introduce a
statistic  on pairs of permutations in $(S_n)^J$ for
$J=S\setminus\{s_{i-2},s_{i-1},s_i\}$. Then we give a formula for
 $R_{u,v}^{J,x}(q)$, where $J=S\setminus\{s_{i-2},s_{i-1},s_i\}$ and
$i$ appears after $i-1$ in $v$.
We also pose a conjecture for
$R_{u,v}^{J,x}(q)$, where $J=S\setminus\{s_{k},s_{k+1},\ldots,s_i\}$ with $1\leq k\leq i\leq n-1$ and the elements $k+1,k+2,\ldots, i$ appear in increasing order in $v$.
\end{abstract}

\noindent{\bf Keywords:} parabolic Kazhdan-Lusztig $R$-polynomial, the symmetric group, Bruhat order

\vspace{6pt}

\noindent{\bf AMS Classifications:} 20F55, 05E99

\section{Introduction}

Parabolic  $R$-polynomials for a Coxeter group   were introduced   by
 Deodhar \cite{Deodhar} as   parabolic analogues
 of ordinary $R$-polynomials defined  by Kazhdan and
 Lusztig \cite{KL}. In this paper, we consider the computation of
parabolic  $R$-polynomials for the symmetric group.
Let $S_n$ be the symmetric group   on $\{1,2,\ldots,n\}$, and let $S=\{s_1,s_2,\ldots, s_{n-1}\}$ be the generating set of $S_n$, where for $1\leq i\leq n-1$, $s_i$
 is the adjacent transposition that interchanges the elements $i$ and $i+1$. For a
  subset  $J\subseteq S$, let $(S_n)_J$ be the  parabolic subgroup
generated by $J$, and let
$(S_n)^{J}$
 be the set of minimal  coset representatives of $S_n/(S_n)_J$.
 Assume that $u$ and $v$ are two permutations in $(S_n)^J$ such that  $u\leq v$ in the Bruhat order. For $x\in \{q,-1\}$,
let $R_{u,v}^{J,x}(q)$ denote the parabolic  $R$-polynomial indexed by $u$ and $v$.
When  $J=S\setminus\{s_i\}$,
Brenti \cite{Brenti-1} found a formula  for  $R_{u,v}^{J,x}(q)$.
Recently, Brenti \cite{Brenti-2}
obtained  an expression for   $R_{u,v}^{J,x}(q)$ for
$J=S\setminus\{s_{i-1}, s_{i}\}$.

In this paper,  we consider the case $J=S\setminus \{s_{i-2},s_{i-1},s_i\}$.
We introduce  a
statistic   on pairs of permutations in $(S_n)^J$
and we give a formula for  $R_{u,v}^{J,x}(q)$, where
$i$ appears after $i-1$ in $v$.
We also  conjecture  a formula  for $R_{u,v}^{J,x}(q)$,
where $J=S\setminus\{s_{k},s_{k+1},\ldots,s_i\}$ with $1\leq k\leq i\leq n-1$ and the elements $k+1,k+2,\ldots, i$ appear in increasing order in $v$.

Let us begin with some terminology and notation.
For a Coxeter group $W$ with a generating set $S$, let
 $T=\{wsw^{-1}\,|\,w\in W,\,s\in S\}$ be
 the set of reflections of $W$. For $w\in W$,
 the length $\ell(w)$ of $w$ is defined
 as  the smallest $k$ such that $w$ can
 be written as a product of $k$ generators in $S$.
For $u,v\in W$, we say that $u\le v$  in
the Bruhat order if  there exists a
sequence $t_1,t_2,\ldots,t_r$ of reflections
 such that $v=u t_1 t_2 \cdots  t_r$ and
 $\ell(u t_1\cdots t_i )>\ell(u  t_1\cdots t_{i-1})$
 for $1\leq i\leq r$.

For a subset  $J\subseteq S$,  let $W_J$ be the parabolic subgroup generated by $J$, and let $W^J$ be the set of minimal right coset representatives of $W/W_J$, that is,
\begin{equation}\label{m-1}
W^J=\{w\in W\,|\,\ell(sw)>\ell(w),\ \text{for all $s\in J$}\}.
\end{equation}
We use $D_{R}(w)$ to denote the set of right descents of $w$, that is,
\begin{equation}\label{m-2}
D_{R}(w)=\{s\in S\,|\, \ell(w s)<\ell(w)\}.
\end{equation}
For  $u, v \in W^J$, the parabolic $R$-polynomial  $R_{u,v}^{J,x}(q)$
can be recursively determined  by the following property.

\begin{theo}[\mdseries{Deodhar \cite{Deodhar}}]\label{cfg-4}
Let $(W,S)$ be a Coxeter system and $J$ be a subset of $S$.
 Then, for each $x\in \{q,-1\}$, there is a unique family $\{R_{u,v}^{J,x}(q)\}_{u,v\in W^J}$ of polynomials
with integer coefficients such that for
all $u,v\in W^J$,
\begin{itemize}
\item[$\mathrm{(i)}$] if $u\nleq v$, then $R_{u,v}^{J,x}(q)=0$;
\item[$\mathrm{(ii)}$] if $u= v$, then $R_{u,v}^{J,x}(q)=1$;
\item[$\mathrm{(iii)}$] if $u<v$, then for any $s\in D_{R}(v)$,
\[R_{u,v}^{J,x}(q)=\left\{
        \begin{array}{lll}
          R_{us,\,vs}^{J,x}(q), & \hbox{\rm{if} $s\in D_{R}(u)$,} \\[5pt]
          qR_{us,\,vs}^{J,x}(q)+(q-1)R_{u,\,vs}^{J,x}(q), & \hbox{\rm{if} $s \notin  D_{R}(u)$ and $us\in W^J$},\\[5pt]
          (q-1-x)R_{u,\,vs}^{J,x}(q), & \hbox{\rm{if} $s \notin  D_{R}(u)$ and $us\not\in W^J$.}
        \end{array}
      \right.
\]
\end{itemize}
\end{theo}

Notice that when $J=\emptyset$, the parabolic $R$-polynomial  $R_{u,v}^{J,x}(q)$ reduces to an ordinary $R$-polynomial  $R_{u,v}(q)$, see, for example, Bj\"{o}rner and Brenti \cite[Chapter 5]{BB} or Humphreys \cite[Chapter 7]{Humphreys}.
The parabolic $R$-polynomials $R_{u,v}^{J,x}(q)$ for $x=q$ and $x=-1$  satisfy the following relation, so that we only need to consider the computation
for the case  $x=q$.

\begin{theo}[\mdseries{Deodhar \cite[Corollary 2.2]{Deodhar-2}}]\label{z}
For  $u,v\in W^J$ with $u\leq v$,
\[q^{\ell(v)-\ell(u)}R_{u,v}^{J,q}\left(\frac{1}{q}\right)=(-1)^{\ell(v)-\ell(u)}R_{u,v}^{J,-1}(q).\]
\end{theo}

There is no  known  explicit formula
for $R_{u,v}^{J,x}(q)$ for  a general Coxeter system $(W,S)$,
 and even for the symmetric group.
When $W=S_n$, Brenti \cite{Brenti-1,Brenti-2}   found
formulas  for  $R_{u,v}^{J,x}(q)$ for certain subsets $J$, namely, $J=S\setminus\{s_i\}$ or $J=S\setminus\{s_{i-1},s_{i}\}$.
To describe the   formulas for the parabolic $R$-polynomials obtained by Brenti \cite{Brenti-1,Brenti-2}, we recall some statistics on pairs of permutations in $(S_n)^J$ with
$J=S\setminus\{s_i\}$ or $J=S\setminus\{s_{i-1},s_{i}\}$.

A permutation $u=u_1u_2\cdots u_n$ in $S_n$ is also considered as
  a bijection
on $\{1,2,\ldots,n\}$ such that $u(i)=u_i$ for $1\leq i\leq n$.
For  $u,v\in S_n$, the product $uv$ of $u$ and $v$ is defined
as a bijection such that $uv(i)=u(v(i))$ for $1\leq i\leq n$.
For $1\leq i\leq n-1$, the adjacent transposition  $s_i$  is the permutation that interchanges the elements $i$ and $i+1$.
The length of a permutation $u\in S_n$ can be interpreted as the number of inversions of $u$, that is,
\begin{equation}\label{len}
\ell(u)=|\{(i,j)\,|\, 1\leq i<j\leq n,\  u(i)>u(j)\}|.
\end{equation}
By \eqref{m-2} and \eqref{len},  the right descent set of a permutation $u\in S_n$ is
given by
\[D_R(u)=\{s_i\,|\, 1\leq i\leq n-1,\ u(i)>u(i+1)\}.\]

When $J=S\setminus\{s_i\}$, it follows from  \eqref{m-1} and \eqref{len} that a permutation $u\in S_n$ belongs to $(S_n)^J$ if and only if   the elements $1,2,\ldots,i$
 as well as the elements $i+1,i+2,\ldots,n$ appear in increasing order in $u$, or equivalently,
\[u^{-1}(1)<u^{-1}(2)<\cdots<u^{-1}(i)\ \ \ \ \text{and}\ \  \ \ u^{-1}(i+1)<u^{-1}(i+2)<\cdots<u^{-1}(n).\]
For $n\geq 1$, we use $[n]$ to denote the set $\{1,2,\ldots,n\}$.
For $J=S\setminus\{s_i\}$ and   $u, v\in (S_n)^J$, let
 \[D(u,v)=v^{-1}([i])\setminus u^{-1}([i]).\]
For $1\leq j\leq n$, let
\begin{equation*}
 a_{j}(u,v)=|\{r\in u^{-1}([i])\,|\, r<j\}|-|\{r\in v^{-1}([i])\,|\, r<j\}|.
 \end{equation*}
Brenti \cite{Brenti-1} obtained the following  formula   for   $R_{u,v}^{J,x}(q)$, where $J=S\setminus\{s_i\}$.

\begin{theo}[\mdseries{Brenti \cite[Corollary 3.2]{Brenti-1}}]\label{cfg-2}
Let $J=S\setminus\{s_i\}$, and let $u, v\in (S_n)^J$ with $u\leq v$. Then
\begin{equation*}
R_{u,v}^{J,q}(q)=(-1)^{\ell(v)-\ell(u)}
\prod_{j\in D(u,v)}\left(1-q^{a_j(u,v)}\right).
\end{equation*}
\end{theo}

We now turn to the case  $J=S\setminus\{s_{i-1},s_i\}$.
In this case,
 it can be seen from \eqref{m-1} and \eqref{len} that a permutation $u\in S_n$ belongs to $(S_n)^J$ if and only if
\[u^{-1}(1)<u^{-1}(2)<\cdots<u^{-1}(i-1)\ \ \ \ \text{and}\ \  \ \ u^{-1}(i+1)<u^{-1}(i+2)<\cdots<u^{-1}(n).\]
For  $u, v\in (S_n)^J$, let
\[\widetilde{D}(u,v)=v^{-1}([i-1])\setminus u^{-1}([i-1]).\]
For  $1\leq j\leq n$, let
\begin{equation*}
\widetilde{ a}_{j}(u,v)=|\{r\in u^{-1}([i-1])\,|\, r<j\}|-|\{r\in v^{-1}([i-1])\,|\, r<j\}|.
 \end{equation*}
The following formula is due to Brenti \cite{Brenti-2}.

\begin{theo}[\mdseries{Brenti \cite[Theorem 3.1]{Brenti-2}}]\label{cfg-3}
Let $J=S\setminus \{s_{i-1}, s_i\}$, and let $u, v\in (S_n)^J$ with $u\leq v$. Then
\begin{equation*}
R_{u,v}^{J,q}(q)=\left\{
        \begin{array}{ll}(-1)^{\ell(v)-\ell(u)}\left(1-q+cq^{1+a_{v^{-1}(i)}(u,v)}\right)
\prod_{j\in D(u,v)}\left(1-q^{a_j(u,v)}\right), & \hbox{\rm{if} $u^{-1}(i)\geq v^{-1}(i)$},\\[7pt]
(-1)^{\ell(v)-\ell(u)}\left(1-q+cq^{1+\widetilde{a}_{v^{-1}(i)}(u,v)}\right)
\prod_{j\in \widetilde{D}(u,v)}\left(1-q^{\widetilde{a}_j(u,v)}\right), & \hbox{\rm{if} $u^{-1}(i)\leq v^{-1}(i)$},
\end{array}
      \right.
\end{equation*}
where $c=\delta_{u^{-1}(i),v^{-1}(i)}$ is the Kronecker delta function.
\end{theo}

It should be noted  that the sets  $(S_n)^J$ for  $J=S\setminus\{s_i\}$ and $J=S\setminus\{s_{i-1},s_{i}\}$ are called
 tight quotients  of $S_n$   by Stembridge \cite{Stembridge} in the study of  the Bruhat order of  Coxeter groups.
Therefore,   combining Theorem  \ref{cfg-2} and Theorem   \ref{cfg-3} leads to an expression for the parabolic
$R$-polynomials for   tight quotients of the symmetric group.

In this paper, we    obtain a formula for
$R_{u,v}^{J,x}(q)$, where
  $J=S\setminus \{s_{i-2},s_{i-1},s_i\}$
  and
$i$ appears after $i-1$ in $v$.
We also conjecture  a more general formula  for $R_{u,v}^{J,x}(q)$,
where $J=S\setminus\{s_{k},s_{k+1},\ldots,s_i\}$ with $1\leq k\leq i\leq n-1$ and the elements $k+1,k+2,\ldots, i$ appear in increasing order in $v$.

\section{A formula for $R_{u,v}^{J,q}(q)$ with $J=S\setminus \{s_{i-2}, s_{i-1}, s_i\}$}

In this section, we present a formula for $R_{u,v}^{J,q}(q)$, where $J=S\setminus \{s_{i-2}, s_{i-1}, s_i\}$ and $v$ is a permutation in $(S_n)^J$ such that $i$ appears after $i-1$ in $v$.
We also give a conjectured  formula  for $R_{u,v}^{J,q}(q)$,
where $J=S\setminus\{s_{k},s_{k+1},\ldots,s_i\}$ with $1\leq k\leq i\leq n-1$ and the elements $k+1,k+2,\ldots, i$ appear in increasing order in $v$.

For $u,v\in (S_n)^J$ with $u\leq v$, our formula for $R_{u,v}^{J,q}(q)$  relies on
a vector of statistics on $(u, v)$, denoted  $(a_1(u,v), a_2(u,v), \ldots, a_n(u,v))$. Notice that a permutation $u\in S_n$ belongs to $(S_n)^J$ if and only if   the elements $1,2,\ldots,i-2$
 as well as the elements $i+1,i+2,\ldots,n$ appear in increasing order in $u$.
To define $a_j(u,v)$, we need to consider the positions of
the elements
$i-1$ and $i$ in  $u$ and $v$.
For convenience, let $u^{-1}=p_1p_2\cdots p_n$ and $v^{-1}=q_1q_2\cdots q_n$,
that is,  $t$ appears in position $p_{t}$ in $u$,
and appears in position  $q_t$ in $v$. The following set $A(u,v)$
is defined based on the relations $p_{i-1}\geq q_{i-1}$ and $p_i\geq q_i$.
More precisely, $A(u,v)$ is a subset of $\{i-1, i\}$ such that $i-1 \in A(u,v)$
if and only if $p_{i-1}\geq q_{i-1}$, and $i\in A(u,v)$ if and only if $p_i\geq q_i$.
Set
\[B(u,v)= \{1,2,\ldots,i-2\} \cup A(u,v).\]

For  $1\leq j\leq n$, we define $a_{j}(u,v)$ to be the number of elements of
$B(u,v)$ that are contained in $\{u_1,\ldots,u_{j-1}\}$ minus the number of
elements of $B(u,v)$ that are contained in   $\{v_1,\ldots,v_{j-1}\}$,
that is,
\begin{equation}\label{cfg-a}
a_{j}(u,v)=|\{r\in u^{-1}(B(u,v))\,|\, r< j\}|-|\{r\in v^{-1}(B(u,v))\,|\, r<j\}|.
\end{equation}
For example, let $n=9$ and  $i=5$, so that $J=S\setminus \{s_3,s_4,s_5\}$.
Let
\begin{equation}\label{la}
u=416273859\ \ \ \text{and}\ \ \ v=671489253
\end{equation}
be two permutations in $(S_9)^J$.
Then we have $A(u,v)=\{5\}$, $B(u,v)=\{1,2,3,5\}$,
and
\begin{equation}\label{lao}
(a_1(u,v),\ldots,a_{9}(u,v))=(0,0,1,0,1,1,2,1,1).
\end{equation}

The following theorem gives a formula for $R_{u,v}^{J,q}(q)$.

\begin{theo}\label{main}
Let $J=S\setminus\{s_{i-2},s_{i-1}, s_i\}$, and let $v$ be a permutation in $(S_n)^J$    such that  $i$ appears after $i-1$ in $v$. Let
\begin{equation}\label{cfg-b}
D(u,v)=v^{-1}(B(u,v))\setminus u^{-1}(B(u,v)).
\end{equation}
Then, for any $u\in (S_n)^J$ with $u\leq v$, we have
\begin{align}\label{cfg-main}
R_{u,v}^{J,q}(q)=(-1)^{\ell(v)-\ell(u)}
&\left(1-q+\delta_{u^{-1}(i-1),v^{-1}(i-1)}q^{1+a_{v^{-1}(i-1)}(u,v)}\right)\nonumber\\[5pt]
&\left(1-q+\delta_{u^{-1}(i),v^{-1}(i)}q^{1+a_{v^{-1}(i)}(u,v)}\right)\prod_{j\in D(u,v)}\left(1-q^{a_j(u,v)}\right).
\end{align}
\end{theo}

Let us give an example for Theorem \ref{main}. Assume that
$u$ and $v$ are the permutations as given in \eqref{la}.
Then we  have $D(u,v)=\{3,7,9\}$. In view of  \eqref{lao},  formula   \eqref{cfg-main}
gives
\[R_{u,v}^{J,q}(q)=(1-q)^3(1-q^2)(1-q+q^2).\]

To prove the above theorem,
 we need  a  criterion for the relation of two  permutations in
$(S_n)^{J}$ with respect to the Bruhat order.
For  $u,v\in (S_n)^{J}$ and $1\leq j\leq n$, let
\begin{align}
b_{1,j}(u,v)=&\left|\{r\in u^{-1}([i-2])\,|\, r<j\}\right|-\left|\{r\in v^{-1}([i-2])\,|\, r<j\}\right|,\label{o-1}\\[5pt]
b_{2,j}(u,v)=&\left|\{r\in u^{-1}([i-1])\,|\, r<j\}\right|-\left|\{r\in v^{-1}([i-1])\,|\, r<j\}\right|,\label{o-2}\\[5pt]
b_{3,j}(u,v)=&\left|\{r\in u^{-1}([i])\,|\, r<j\}\right|-\left|\{r\in v^{-1}([i])\,|\, r<j\}\right|.\label{o-3}
\end{align}
The following proposition shows that we can use $b_{h,j}(u,v)$
with  $h=1,2,3$ and $1\leq j\leq n$ to determine
whether $u\leq v$ in the Bruhat order.

\begin{prop}\label{prop-2}
Let $J=S\setminus \{s_{i-2}, s_{i-1},s_i\}$, and let $u,v\in (S_n)^{J}$.
Then,
$u\leq v$ if and only if $b_{h,j}(u,v)\geq 0$ for $h=1,2,3$ and $1\leq j\leq n$.
\end{prop}

\pf Recall that for $u,v\in S_n$,   the following   statements are equivalent \cite[Chapter 2]{BB}:
\begin{itemize}
\item[(i)] $u\leq v$;
\item[(ii)] $u^{-1}\leq v^{-1}$;
\item[(iii)] for any $1\leq j\leq n$ and any $1\leq t\leq n$,
\[\left|\{r\in u([t])\,|\, r<j\}\right|\geq \left|\{r\in v([t])\,|\, r<j\}\right|.\]
\item[(iv)] for any $1\leq j\leq n$ and any $t$ such that $s_t\in D_R(u)$,
\[\left|\{r\in u([t])\,|\, r<j\}\right|\geq \left|\{r\in v([t])\,|\, r<j\}\right|.\]
\end{itemize}
Therefore, it suffices to show that for $u,v\in (S_n)^{J}$,
$u^{-1}\leq v^{-1}$ if and only if $b_{h,j}(u,v)\geq 0$ for
  $h=1,2,3$ and $1\leq j\leq n$.

Notice  that the elements $1,2,\ldots,i-2$
 as well as the elements $i+1,i+2,\ldots,n$ appear in increasing order in $u$.
We see that if  $s_t\in D_R(u^{-1})$, then $t$ equals $i-2$, $i-1$, or $i$. Thus,  $u^{-1}\leq v^{-1}$
if and only if  for $h=1,2,3$ and $1\leq j\leq n$,
 \[\left|\{r\in u^{-1}([i-3+h])\,|\, r<j\}\right|\geq \left|\{r\in v^{-1}([i-3+h])\,|\, r<j\}\right|,\]
or equivalently, $u^{-1}\leq v^{-1}$ if and only if $b_{h,j}(u,v)\geq 0$ for $h=1,2,3$ and $1\leq j\leq n$.
 This completes the proof.
\qed

We are now in a position to present a proof of Theorem \ref{main}.

\noindent
{\it {Proof of Theorem \ref{main}.}}
Assume that $J=S\setminus \{s_{i-2},s_{i-1},s_i\}$, and $u$ and $v$ are two permutations
in $(S_n)^J$ such that $u\leq v$.
Write $u^{-1}=p_1p_2\cdots p_n$ and $v^{-1}=q_1q_2\cdots q_n$.
By the definitions of  $(a_1(u,v),\ldots, a_n(u,v))$ and $D(u,v)$,
 we  consider   the following four cases:
\begin{align}
&p_{i-1}\geq q_{i-1}\ \ \ \ \text{and}\ \ \ \ p_{i}\geq q_{i},\label{pr-a}\\[5pt]
&p_{i-1}\geq q_{i-1}\ \ \ \ \text{and}\ \ \ \ p_{i}< q_{i},\label{pr}\\[5pt]
&p_{i-1}< q_{i-1}\ \ \ \ \text{and}\ \ \ \ p_{i}\geq q_{i},\label{pr-b}\\[5pt]
&p_{i-1}< q_{i-1}\ \ \ \ \text{and}\ \ \ \ p_{i}< q_{i}.\label{pr-c}
\end{align}

We conduct induction on  $\ell(v)$. When $\ell(v)=0$, formula  \eqref{cfg-main} is
obvious. Assume that $\ell(v)>0$ and formula \eqref{cfg-main} is true for $\ell(v)-1$.
We proceed to prove   \eqref{cfg-main} for $\ell(v)$.
We shall only provide a proof for the case   in \eqref{pr}. The other cases can be
justified by using similar arguments.
By \eqref{cfg-a} and \eqref{pr}, we see that for $1\leq k\leq n$,
\begin{equation}\label{cfg-24}
a_k(u,v)=|\{r\in u^{-1}([i-1])\,|\, r<k\}|-|\{r\in v^{-1}([i-1])\,|\, r<k\}|.
\end{equation}
Moreover, by \eqref{cfg-b} and \eqref{pr} we find that
\begin{equation}\label{cfg-25}
D(u,v)=v^{-1}([i-1])\setminus  u^{-1}([i-1]).
\end{equation}

Let  $s_j\in D_R(v)$ be a right descent of $v$,
 that is, $v(j)>v(j+1)$, where $1\leq j\leq n-1$.
Keep in mind that $i$ appears after $i-1$ in $v$, namely, $q_i>q_{i-1}$, and that
the elements $1,2,\ldots,i-2$
 as well as the elements $i+1,i+2,\ldots,n$ appear in increasing order in $v$.
So we get all possible choices  of $v(j)$ and $v(j+1)$   as listed in Table \ref{ttt}.
\begin{table}[h,t]
\begin{center}
\begin{tabular}{|l|l|l|l|}\hline
$1$&$v(j)>i$\ \  and\ \  $v(j+1)=i$\\\hline
$2$&$v(j)>i$\ \  and\ \  $v(j+1)=i-1$\\\hline
$3$&$v(j)>i$\ \  and\ \  $v(j+1)<i-1$\\\hline
$4$&$v(j)=i$\ \  and\ \  $v(j+1)<i-1$\\\hline
$5$&$v(j)=i-1$\ \  and\ \  $v(j+1)<i-1$\\\hline
\end{tabular}
\end{center}
\caption{The choices  of $v(j)$ and $v(j+1)$ in $v$.}\label{ttt}
\end{table}

According to whether $s_j$ is a right descent of $u$,
we have the following two cases.

\noindent
Case 1:  $s_j\in D_R(u)$, that is, $u(j)>u(j+1)$.
Since the elements $1,2,\ldots,i-2$
 as well as the elements $i+1,i+2,\ldots,n$ appear in increasing order in $u$,
the possible choices of $u(j)$ and $u(j+1)$ are as given in
 Table \ref{tt1}.
\begin{table}[h,t]
\begin{center}
\begin{tabular}{|l|l|l|l|}\hline
1&$u(j)>i$\ \  and\ \ $u(j+1)=i$\\\hline
2&$u(j)>i$\ \  and\ \  $u(j+1)=i-1$\\\hline
3&$u(j)>i$\ \  and\ \  $u(j+1)<i-1$\\\hline
4&$u(j)=i$\ \  and\ \  $u(j+1)=i-1$\\\hline
5&$u(j)=i$\ \  and\ \  $u(j+1)<i-1$\\\hline
6&$u(j)=i-1$\ \  and\ \  $u(j+1)<i-1$\\\hline
\end{tabular}
\end{center}
\caption{The choices  of $u(j)$ and $u(j+1)$ in $u$ in Case 1.}\label{tt1}
\end{table}

We only give a proof for the case when $v$ satisfies
 Condition 5 in Table \ref{ttt}, that is,
$v(j)=i-1$ and $v(j+1)<i-1$. The remaining cases can be dealt with  in
the same manner. By the classification of $u$ in Table \ref{tt1}, there are six
 subcases to consider.

\noindent
Subcase 1: $u(j)>i$ and $u(j+1)=i$.
By \eqref{cfg-a}, it is easy to check that for $1\leq k\leq n$,
\[a_k(us_j,vs_j)=a_k(u,v).\]
Moreover, it follows from \eqref{cfg-b} that
\[D(us_j,vs_j)=D(u,v).\]
By the induction hypothesis, we deduce that
\begin{align*}
R_{u,v}^{J,q}(q)&=R_{us_j,vs_j}^{J,q}(q)\\[5pt]
&=(-1)^{\ell(vs_j)-\ell(us_j)} (1-q)^2
\prod_{k\in D(us_j,vs_j)}\left(1-q^{a_k(us_j,vs_j)}\right)\\[5pt]
&=(-1)^{\ell(v)-\ell(u)} (1-q)^2
\prod_{k\in D(u,v)}\left(1-q^{a_k(u,v)}\right),
\end{align*}
as required.

\noindent
Subcase 2: $u(j)>i$ and $u(j+1)=i-1$.
Notice that in this case $us_j$ and $vs_j$  satisfy the relation in \eqref{pr-c}.
So we have $B(us_j,vs_j)=[i-2]$.
By \eqref{cfg-a} and \eqref{cfg-b}, it is easily verified  that for $1\leq k\leq n$,
\[a_k(us_j,vs_j)=a_k(u,v),\]
and
\begin{align*}
D(us_j,vs_j)=& (vs_j)^{-1}([i-2])\setminus (us_j)^{-1}([i-2])\\[5pt]
=& v^{-1}([i-1])\setminus u^{-1}([i-1])\\[5pt]
=&D(u,v).
\end{align*}
By  the induction hypothesis, we get
\begin{align*}
R_{u,v}^{J,q}(q)=&R_{us_j,vs_j}^{J,q}(q)\\[5pt]
=&(-1)^{\ell(vs_j)-\ell(us_j)} (1-q)^2
 \prod_{k\in D(us_j,vs_j)}\left(1-q^{a_k(us_j,vs_j)}\right)\\[5pt]
=&(-1)^{\ell(v)-\ell(u)} (1-q)^2
\prod_{k\in D(u,v)}\left(1-q^{a_k(u,v)}\right),
\end{align*}
which implies     \eqref{cfg-main}.

\noindent
Subcase 3: $u(j)>i$ and $u(j+1)<i-1$.
We find that
\[
a_{j+1}(us_j,vs_j)=a_{j+1}(u,v)+1\ \ \text{and}\ \
a_k(us_j,vs_j)=a_k(u,v),\ \ \text{for}\ k\neq j+1,
\]
and
\[D(us_j,vs_j)=\left(D(u,v)\setminus \{j\}\right)\cup \{j+1\}.\]
Thus, the induction hypothesis yields that
\begin{align*}
R_{u,v}^{J,q}(q)=&R_{us_j,vs_j}^{J,q}(q)\\[5pt]
=&(-1)^{\ell(vs_j)-\ell(us_j)} (1-q)^2
 \prod_{k\in D(us_j,vs_j)}\left(1-q^{a_k(us_j,vs_j)}\right)\\[5pt]
=&(-1)^{\ell(v)-\ell(u)} (1-q)^2
\frac{1-q^{a_{j+1}(us_j,vs_j)}}{1-q^{a_{j}(u,v)}}\prod_{k\in D(u,v)}\left(1-q^{a_k(u,v)}\right),
\end{align*}
which
reduces to    \eqref{cfg-main}, since
\[a_{j+1}(us_j,vs_j)=a_{j}(u,v).\]

\noindent
Subcase 4: $u(j)=i$ and $u(j+1)=i-1$. This case is similar to
   Subcase 2, and hence the proof is omitted.

\noindent
Subcase 5: $u(j)=i$ and $u(j+1)<i-1$. This case is similar to Subcase 3.

\noindent
Subcase  6: $u(j)=i-1$ and $u(j+1)<i-1$.
For $1\leq k\leq n$, we have
\[
a_k(us_j,vs_j)=a_k(u,v)
\]
and
\[D(us_j,vs_j)=D(u,v).\]
By the induction hypothesis,  we find that
\begin{align}
R_{u,v}^{J,q}(q)=&R_{us_j,vs_j}^{J,q}(q)\nonumber\\[5pt]
=&(-1)^{\ell(vs_j)-\ell(us_j)} \left(1-q+q^{1+a_{j+1}(us_j,vs_j)}\right)(1-q)
 \prod_{k\in D(us_j,vs_j)}\left(1-q^{a_k(us_j,vs_j)}\right).\label{yx}
\end{align}
Noticing the following relation
\[a_{j+1}(us_j,vs_j)=a_{j}(u,v),\]
formula \eqref{yx} can be rewritten as
\[
R_{u,v}^{J,q}(q)=(-1)^{\ell(v)-\ell(u)} \left(1-q+q^{1+a_{j}(u,v)}\right)
\left(1-q\right) \prod_{k\in D(u,v)}\left(1-q^{a_k(u,v)}\right),
\]
as required.

\vspace{6pt}

\noindent
Case 2:  $s_j\not\in D_R(u)$, that is, $u(j)<u(j+1)$.
The possible choices of $u(j)$ and $u(j+1)$ are given in
 Table \ref{tt2}.
\begin{table}[h,t]
\begin{center}
\begin{tabular}{|l|l|l|l|}\hline
1&$u(j)<u(j+1)<i-1$\\\hline
2&$u(j)<u(j+1)=i-1$\\\hline
3&$u(j)<i-1$\ \  and\ \  $u(j+1)=i$\\\hline
4&$u(j)=i-1$\ \  and\ \  $u(j+1)=i$\\\hline
5&$u(j)<i-1$\ \  and\ \  $u(j+1)>i$\\\hline
6&$u(j)=i-1$\ \  and\ \  $u(j+1)>i$\\\hline
7&$i=u(j)<u(j+1)$\\\hline
8&$i<u(j)<u(j+1)$\\\hline
\end{tabular}
\end{center}
\caption{The choices  of $u(j)$ and $u(j+1)$ in $u$ in Case 2.}\label{tt2}
\end{table}

We shall provide proofs for three subcases: (i) $v$ satisfies Condition 1 in Table \ref{ttt} and $u$ satisfies Condition 7 in Table \ref{tt2};
(ii) $v$ satisfies Condition 3 in Table \ref{ttt} and $u$ satisfies Condition 3 in Table \ref{tt2}; (iii) $v$ satisfies Condition 5 in Table \ref{ttt} and $u$ satisfies Condition 3 in Table \ref{tt2}.
The verifications in other situations are similar or relatively easier.

\noindent
Subcase (i):  $v(j)>i$, $v(j+1)=i$, $i=u(j)<u(j+1)$.
By Theorem \ref{cfg-4}, we have
\begin{equation}\label{qq-1}
R_{u,v}^{J,q}(q)=qR_{us_j,vs_j}^{J,q}(q)+(q-1)R_{u,vs_j}^{J,q}(q).
\end{equation}
We need to compute   $R_{us_j,vs_j}^{J,q}(q)$ and $R_{u,vs_j}^{J,q}(q)$.
We first compute $R_{u,vs_j}^{J,q}(q)$.
Notice that $u$ and $vs_j$  satisfy the relation in \eqref{pr-a}.
By \eqref{cfg-a},
we obtain  that for $1\leq k\leq n$,
\begin{align}
a_{k}(u,vs_j)=&\left|\{r\in u^{-1}([i])\,|\, r<k\}\right|
-\left|\{r\in (vs_j)^{-1}([i])\,|\, r<k\}\right|\nonumber\\[5pt]
=&\left|\{r\in u^{-1}([i-1])\,|\, r<k\}\right|
-\left|\{r\in v^{-1}([i-1])\,|\, r<k\}\right|\nonumber\\[5pt]
=&a_{k}(u,v).\nonumber
\end{align}
Moreover, by \eqref{cfg-b} we have
\begin{align*}
D(u,vs_j)&=(vs_j)^{-1}([i])\setminus u^{-1}([i])\\[5pt]
&=v^{-1}([i-1])\setminus u^{-1}([i-1])\\[5pt]
&=D(u,v).
\end{align*}
By the induction hypothesis, we deduce that
\begin{align}
R_{u,vs_j}^{J,q}(q)=& (-1)^{\ell(vs_j)-\ell(u)}
\left(1-q+\delta_{u^{-1}(i-1),(vs_j)^{-1}(i-1)}q^{1+a_{(vs_j)^{-1}(i-1)}
(u,vs_j)}\right)\nonumber\\[5pt]
&\ \ \left(1-q+q^{1+a_{j}(u,vs_j)}\right)\prod_{k\in D(u,vs_j)}\left(1-q^{a_k(u,vs_j)}\right)\nonumber\\[5pt]
=& (-1)^{\ell(v)-\ell(u)-1}
\left(1-q+\delta_{u^{-1}(i-1),v^{-1}(i-1)}q^{1+a_{v^{-1}(i-1)}(u,v)}\right)\nonumber\\[5pt]
&\ \ \left(1-q+q^{1+a_{j}(u,v)}\right)\prod_{k\in D(u,v)}\left(1-q^{a_k(u,v)}\right). \label{cfg-336}
\end{align}

To compute  $R_{us_j,vs_j}^{J,q}(q)$, we consider two cases according to  whether $us_j\leq vs_j$.  First, we assume that $us_j\leq vs_j$.
Since  $us_j$ and $vs_j$  satisfy the relation in \eqref{pr-a}, by \eqref{cfg-a} we see that
\begin{equation} \label{y-2}
a_{j+1}(us_j,vs_j)=a_{j+1}(u,v)-1 \ \ \text{and }\ \  a_k(us_j,vs_j)=a_k(u,v),\ \ \text{for $k\neq j+1$} .
\end{equation}
Moreover, by \eqref{cfg-b} we get
\begin{align}
D(us_j,vs_j)&=(vs_j)^{-1}([i])\setminus (us_j)^{-1}([i])\nonumber\\[5pt]
&=D(u,v)\cup \{j\}.\label{y-1}
\end{align}
Combining \eqref{y-2} and  \eqref{y-1} and  applying  the induction hypothesis,
we deduce that
\begin{align}
R_{us_j,vs_j}^{J,q}(q)=& (-1)^{\ell(vs_j)-\ell(us_j)}
\left(1-q+\delta_{(us_j)^{-1}(i-1),(vs_j)^{-1}(i-1)}q^{1+a_{(vs_j)^{-1}(i-1)}
(us_j,vs_j)}\right)\nonumber\\[5pt]
&\ \ (1-q)\prod_{k\in D(us_j,vs_j)}\left(1-q^{a_k(us_j,vs_j)}\right)\nonumber\\[5pt]
=& (-1)^{\ell(v)-\ell(u)}
\left(1-q+\delta_{u^{-1}(i-1),v^{-1}(i-1)}q^{1+a_{v^{-1}(i-1)}(u,v)}\right)\nonumber\\[5pt]
&\ \ (1-q)\left(1-q^{a_j(u,v)}\right)\prod_{k\in D(u,v)}\left(1-q^{a_k(u,v)}\right).\label{cfg-337}
\end{align}
Substituting \eqref{cfg-336} and \eqref{cfg-337} into \eqref{qq-1}, we obtain that
\begin{align}
R_{u,v}^{J,q}(q)=&qR_{us_j,vs_j}^{J,q}(q)+(q-1)R_{u,vs_j}^{J,q}(q)\nonumber\\[5pt]
=& (-1)^{\ell(v)-\ell(u)}\left(q\left(1-q^{a_j(u,v)}\right)+\left(1-q+q^{1+a_{j}(u,v)}\right)\right)
\nonumber\\[5pt]
&\ \ (1-q)
\left(1-q+\delta_{u^{-1}(i-1),v^{-1}(i-1)}q^{1+a_{v^{-1}(i-1)}(u,v)}\right)\prod_{k\in D(u,v)}\left(1-q^{a_k(u,v)}\right)\nonumber\\[5pt]
=&(-1)^{\ell(v)-\ell(u)}(1-q)
\left(1-q+\delta_{u^{-1}(i-1),v^{-1}(i-1)}q^{1+a_{v^{-1}(i-1)}(u,v)}\right)\nonumber\\[5pt]
&\ \ \prod_{k\in D(u,v)}\left(1-q^{a_k(u,v)}\right).\nonumber
\end{align}

We now consider the case $us_j\not\leq vs_j$.
In this case, we claim that
\begin{equation}\label{at}
a_{j}(u,v)=0.
\end{equation}
 By \eqref{o-1}--\eqref{o-3}, it can be checked that for $1\leq k\leq n$,
\[
b_{1, k}(us_j,vs_j)=b_{1, k}(u,v)\ \ \text{and}\ \ b_{2, k}(us_j,vs_j)=b_{2, k}(u,v),
\]
and
\[b_{3, j+1}(us_j,vs_j)=b_{3, j+1}(u,v)-2\ \   \text{and}\ \ b_{3, k}(us_j,vs_j)=b_{3, k}(u,v),\ \ \text{for}\ k\neq j+1.\]
Since  $us_j\not\leq vs_j$, by Proposition \ref{prop-2}, we see that $b_{3, j+1}(u,v)-2<0$.
On the other hand, since  $j+1\in v^{-1}([i])$ but  $j+1\not\in u^{-1}([i])$,
we  have $b_{3, j+1}(u,v)>0$. So we get $b_{3, j+1}(u,v)=1$.
Therefore,
\[a_{j}(u,v)=b_{3, j+1}(u,v)-1=0.\]
This proves the claim in \eqref{at}.

Combining  \eqref{cfg-336} and \eqref{at}, we obtain that
\begin{align*}
R_{u,v}^{J,q}(q)=&(q-1)R_{u,vs_j}^{J,q}(q) \\[5pt]
=& (-1)^{\ell(v)-\ell(u)}(1-q)
\left(1-q+\delta_{u^{-1}(i-1),v^{-1}(i-1)}q^{1+a_{v^{-1}(i-1)}(u,v)}\right)
\prod_{k\in D(u,v)}\left(1-q^{a_k(u,v)}\right).
\end{align*}

\noindent
Subcase (ii): $v(j)>i$, $v(j+1)<i-1$,  $u(j)<i-1$ and $u(j+1)=i$.
By Theorem \ref{cfg-4}, we have
\begin{equation}\label{mm-1}
R_{u,v}^{J,q}(q)=qR_{us_j,vs_j}^{J,q}(q)+(q-1)R_{u,vs_j}^{J,q}(q).
\end{equation}
We need to compute   $R_{us_j,vs_j}^{J,q}(q)$
and $R_{u,vs_j}^{J,q}(q)$. We first compute $R_{u,vs_j}^{J,q}(q)$.
Using \eqref{cfg-a}, we get
\[
a_{j+1}(u,vs_j)=a_{j+1}(u,v)-1 \ \ \text{and }\ \  a_k(u,vs_j)=a_k(u,v),\ \ \text{for $k\neq j+1$} .
\]
Moreover, by \eqref{cfg-b} we have
\[
D(u,vs_j)=D(u,v)\setminus \{j+1\}.
\]
By  the induction hypothesis,
we deduce that
\begin{align}
R_{u,vs_j}^{J,q}(q)=& (-1)^{\ell(vs_j)-\ell(u)}
\left(1-q+\delta_{u^{-1}(i-1),(vs_j)^{-1}(i-1)}q^{1+a_{(vs_j)^{-1}(i-1)}
(u,vs_j)}\right)
\nonumber\\[5pt]
&\ \ (1-q)\prod_{k\in D(u,vs_j)}\left(1-q^{a_k(u,vs_j)}\right)\nonumber\\[5pt]
=& (-1)^{\ell(v)-\ell(u)-1}
\left(1-q+\delta_{u^{-1}(i-1),v^{-1}(i-1)}q^{1+a_{v^{-1}(i-1)}
(u,v)}\right)(1-q)
\nonumber\\[5pt]
&\ \ \frac{1}{1-q^{a_{j+1}(u,v)}}\prod_{k\in D(u,v)}\left(1-q^{a_k(u,v)}\right). \label{aaa-1}
\end{align}

To compute  $R_{us_j,vs_j}^{J,q}(q)$, we consider two cases according to  whether $us_j\leq vs_j$.  First, we assume that $us_j\leq vs_j$.
In view of \eqref{cfg-a}, it is easy to check that
\[
a_{j+1}(us_j,vs_j)=a_{j+1}(u,v)-2 \ \ \text{and }\ \  a_k(us_j,vs_j)=a_k(u,v),\ \ \text{for $k\neq j+1$}.
\]
Moreover, it follows from \eqref{cfg-b}  that
\[
D(us_j,vs_j)=\left(D(u,v)\setminus \{j+1\}\right)\cup \{j\}.
\]
By  the induction hypothesis,
we obtain that
\begin{align}
R_{us_j,vs_j}^{J,q}(q)=& (-1)^{\ell(vs_j)-\ell(us_j)}
\left(1-q+\delta_{(us_j)^{-1}(i-1),(vs_j)^{-1}(i-1)}q^{1+a_{(vs_j)^{-1}(i-1)}
(us_j,vs_j)}\right)
\nonumber\\[5pt]
&\ \ (1-q)\prod_{k\in D(us_j,vs_j)}\left(1-q^{a_k(us_j,vs_j)}\right)\nonumber\\[5pt]
=& (-1)^{\ell(v)-\ell(u)}
\left(1-q+\delta_{u^{-1}(i-1),v^{-1}(i-1)}q^{1+a_{v^{-1}(i-1)}
(u,v)}\right)(1-q)
\nonumber\\[5pt]
&\ \ \frac{1-q^{a_{j}(us_j,vs_j)}}{1-q^{a_{j+1}(u,v)}}\prod_{k\in D(u,v)}\left(1-q^{a_k(u,v)}\right). \label{aaa-2}
\end{align}
Substituting \eqref{aaa-1} and \eqref{aaa-2} into \eqref{mm-1} and noticing the following relation
\[a_{j}(us_j,vs_j)=a_{j+1}(u,v)-1,\]
we are led to  formula   \eqref{cfg-main}.

We now consider the case $us_j\not\leq vs_j$.
In this case, we  claim that
\begin{equation}\label{aaa-4}
a_{j+1}(u,v)=1.
\end{equation}
By \eqref{o-1}--\eqref{o-3}, it is easily seen that
\begin{align}
&b_{1, j+1}(us_j,vs_j)=b_{1, j+1}(u,v)-2\ \   \text{and}\ \ b_{1, k}(us_j,vs_j)=b_{1, k}(u,v),\ \ \text{for}\ k\neq j+1,\label{aaa-8}\\[5pt]
&b_{2, j+1}(us_j,vs_j)=b_{2, j+1}(u,v)-2\ \   \text{and}\ \ b_{2, k}(us_j,vs_j)=b_{2, k}(u,v),\ \ \text{for}\ k\neq j+1,\\[5pt]
&b_{3, j+1}(us_j,vs_j)=b_{3, j+1}(u,v)-1\ \   \text{and}\ \ b_{3, k}(us_j,vs_j)=b_{3, k}(u,v),\ \ \text{for}\ k\neq j+1.\label{aaa-9}
\end{align}
It is clear that $a_{j+1}(u,v)=b_{2, j+1}(u,v)$. So   the claim in \eqref{aaa-4} reduces to
\[
b_{2, j+1}(u,v)=1.
\]

Since  $j\not\in v^{-1}([i-1])$ but  $j \in u^{-1}([i-1])$,
we have $b_{2, j+1}(u,v)>0$.  Suppose to the contrary that $b_{2, j+1}(u,v)>1$.
In the notation  $u^{-1}=p_1p_2\cdots p_n$ and $v^{-1}=q_1q_2\cdots q_n$, we have  the following two cases.

\noindent
Case (a): $p_{i-1}<j$. By \eqref{pr}, we see that $q_{i-1}<j$ and 
\[b_{1, j+1}(u,v)=b_{2, j+1}(u,v)>1.\]
On the other hand, since $j\not\in v^{-1}([i])$ but   $j \in u^{-1}([i])$, we have $b_{3, j+1}(u,v)>0$. Hence  we conclude that $b_{h,k}(us_j,vs_j)\geq 0$   for
$h=1,2,3$ and $1\leq k\leq  n$. By Proposition \ref{prop-2}, we get $us_j\leq vs_j$, contradicting the assumption $us_j\not\leq vs_j$.

\noindent
Case (b): $p_{i-1}>j$. In this case, 
we find that 
if $q_{i-1}>j$, then 
\[b_{1, j+1}(u,v)= b_{2, j+1}(u,v)>1,\]
whereas if $q_{i-1}<j$, then
\[b_{1, j+1}(u,v)>b_{2, j+1}(u,v)>1.\]
Note that in Case (a), we have shown that  $b_{3, j+1}(u,v)>0$.
So, we obtain that $b_{h,k}(us_j,vs_j)\geq 0$   for
$h=1,2,3$ and $1\leq k\leq  n$. Thus we have   $us_j\leq vs_j$,
contradicting the assumption $us_j\not\leq vs_j$.
This proves the claim in \eqref{aaa-4}.
Substituting \eqref{aaa-4} into   \eqref{aaa-1},
 we arrive at \eqref{cfg-main}.

\vspace{6pt}
\noindent
Subcase (iii): $v(j)=i-1$, $v(j+1)<i-1$, $u(j)<i-1$ and   $u(j+1)=i$.
By Theorem \ref{cfg-4}, we have
\begin{equation}\label{cfg-23}
R_{u,v}^{J,q}(q)=qR_{us_j,vs_j}^{J,q}(q)+(q-1)R_{u,vs_j}^{J,q}(q).
\end{equation}
We need to compute
$R_{us_j,vs_j}^{J,q}(q)$ and $R_{u,vs_j}^{J,q}(q)$.
By \eqref{cfg-a}, we see that for $1\leq k\leq n$,
\[
a_k(u,vs_j)=a_k(u,v).
\]
Moreover, by \eqref{cfg-b} we have
\[
D(u,vs_j)=D(u,v).
\]
By the induction  hypothesis, we obtain that
\begin{align}
R_{u,vs_j}^{J,q}(q)=& (-1)^{\ell(vs_j)-\ell(u)}
(1-q)^2\prod_{k\in D(u,vs_j)}\left(1-q^{a_k(u,vs_j)}\right)\nonumber\\[5pt]
=& (-1)^{\ell(v)-\ell(u)-1}
(1-q)^2\prod_{k\in D(u,v)}\left(1-q^{a_k(u,v)}\right).\label{cfg-27}
\end{align}

To  compute  $R_{us_j,vs_j}^{J,q}(q)$, we  claim  that $us_j\leq vs_j$.
By \eqref{o-1}--\eqref{o-3}, we see that
\begin{align}
&b_{1, j+1}(us_j,vs_j)=b_{1, j+1}(u,v)-2\ \   \text{and}\ \ b_{1, k}(us_j,vs_j)=b_{1, k}(u,v),\ \ \text{for}\ k\neq j+1,\label{z-1}\\[5pt]
&b_{2, j+1}(us_j,vs_j)=b_{2, j+1}(u,v)-1\ \   \text{and}\ \ b_{2, k}(us_j,vs_j)=b_{2, k}(u,v),\ \ \text{for}\ k\neq j+1,\label{z-2}\\[5pt]
&b_{3, k}(us_j,vs_j)=b_{3, k}(u,v),\ \ \text{for $1\leq k\leq n$}.\label{z-3}
\end{align}
Since $j+1\in v^{-1}([i-1])$ but  $j+1\not\in u^{-1}([i-1])$, we   have $b_{2, j+1}(u,v)>0$, which implies that
\begin{equation}\label{z-4}
b_{2, j+1}(us_j,vs_j)=b_{2, j+1}(u,v)-1\geq 0.
\end{equation}
Moreover, since
$p_{i-1}\geq q_{i-1}=j$, we have $p_{i-1}>j$.
So, we deduce that
\[
b_{1, j+1}(u,v)=b_{2, j+1}(u,v)+1>1,
\]
and hence
\begin{equation}\label{z-5}
b_{1, j+1}(us_j,vs_j)=b_{1, j+1}(u,v)-2\geq 0.
\end{equation}
Therefore, for   $h=1,2,3$ and  $1\leq k\leq n$,
\[b_{h, k}(us_j,vs_j)\geq 0,\]
which together with Proposition \ref{prop-2} yields that $us_j\leq vs_j$.
This proves the claim.

By \eqref{cfg-a} and \eqref{cfg-b}, it is easily verified that
\[a_{j+1}(us_j,vs_j)=a_{j+1}(u,v)-1 \ \ \text{and }\ \
a_k(us_j,vs_j)=a_k(u,v),\ \ \text{for $k\neq j+1$}
\]
and
\[
D(us_j,vs_j)=\left(D(u,v)\setminus \{j+1\}\right)\cup\{j\}.
\]
By the induction hypothesis, we deduce that
\begin{align}
R_{us_j,vs_j}^{J,q}(q)=& (-1)^{\ell(vs_j)-\ell(us_j)}
(1-q)^2\prod_{k\in D(us_j,vs_j)}\left(1-q^{a_k(us_j,vs_j)}\right)\nonumber\\[5pt]
=& (-1)^{\ell(v)-\ell(u)}
(1-q)^2\frac{1-q^{a_j(us_j,vs_j)}}{1-q^{a_{j+1}(u,v)}}\prod_{k\in D(u,v)}\left(1-q^{a_k(u,v)}\right).\label{cfg-26}
\end{align}
Since  $a_j(us_j,vs_j)=a_{j+1}(u,v)$, formula   \eqref{cfg-26} becomes
\begin{equation}\label{cfg-28}
R_{us_j,vs_j}^{J,q}(q)=(-1)^{\ell(v)-\ell(u)}
(1-q)^2
\prod_{k\in D(u,v)}\left(1-q^{a_k(u,v)}\right)
\end{equation}
Substituting \eqref{cfg-27} and \eqref{cfg-28} into \eqref{cfg-23},
we are led to  \eqref{cfg-main}.
This completes the proof.
\qed

We conclude this paper with a conjectured  formula for  $R_{u,v}^{J,q}(q)$, where
\[J=S\setminus \{s_{k}, s_{k+1},\ldots,s_i\}\]
with  $1\leq k\leq  i\leq n-1$ and   $k+1,k+2,\ldots,i$
appear in increasing order in $v$.
By \eqref{m-1} and \eqref{len}, a permutation $u\in S_n$ belongs to $(S_n)^{J}$ if and only if the elements $1,2,\ldots,k$
 as well as the elements $i+1,i+2,\ldots,n$ appear in increasing order in $u$.
Let $u,v$ be two permutations in $(S_n)^{J}$.
Write $u^{-1}=p_1p_2\cdots p_n$ and $v^{-1}=q_1q_2\cdots q_n$.
Let
\[A(u,v)=\{t\,|\,k+1\leq t\leq i,\,p_{t}\geq q_{t} \}.\]
Set $B(u,v)$ to
 be the union of $\{1,2,\ldots,k\}$ and $A(u,v)$.
Based on the set $B(u,v)$,
we define $a_{j}(u,v)$ and $D(u,v)$ in the same way as
in \eqref{cfg-a} and \eqref{cfg-b},
respectively.

The following conjecture  has been verified for $n\leq 8$.

\begin{conj}\label{conj}
Let $J=S\setminus \{s_{k}, s_{k+1},\ldots,s_i\}$, and  $v$ be a permutation in  $(S_n)^{J}$  such that $k+1,k+2,\ldots,i$
appear in increasing order in $v$.
Then, for any $u\in (S_n)^{J}$ with $u\leq v$, we have
\begin{equation*}\label{cfg-11}
R_{u,v}^{J,q}(q)=(-1)^{\ell(v)-\ell(u)}
\prod_{t=k+1}^{i}\left(1-q+\delta_{u^{-1}(t),v^{-1}(t)}q^{1+a_{v^{-1}(t)}(u,v)}\right)
\prod_{j\in D(u,v)}\left(1-q^{a_j(u,v)}\right).
\end{equation*}
\end{conj}

  Conjecture \ref{conj} contains  Theorems \ref{cfg-2},  \ref{cfg-3} and
\ref{main} as special cases. When $i=n-1$ and $k=1$, Conjecture \ref{conj} becomes a
 conjectured formula for ordinary $R$-polynomials $R_{u,v}(q)$, that is, if $v$ is a permutation in $S_n$ such that $2,3,\ldots,n-1$
appear in increasing order in $v$,
then for any  $u\leq v $,
\[
R_{u,v}(q)=(-1)^{\ell(v)-\ell(u)}
\prod_{t=2}^{n-1}\left(1-q+\delta_{u^{-1}(t),v^{-1}(t)}q^{1+a_{v^{-1}(t)}(u,v)}\right)
\prod_{j\in D(u,v)}\left(1-q^{a_j(u,v)}\right).
\]

\vspace{.2cm} \noindent{\bf Acknowledgments.}
This work was
supported by the 973 Project, the PCSIRT Project of the Ministry of
Education, and the National Science Foundation of China.


\begin{thebibliography}{99}

\bibitem{BB}
A. Bj\"{o}rner and F. Brenti, Combinatorics of Coxeter groups,
Graduate Texts in Mathematics, Vol. 231, Springer-Verlag, New York,
2005.

\bibitem{Brenti-1}
F. Brenti, Kazhdan-Lusztig and $R$-polynomials, Young's lattice, and Dyck partitions,
Pacific J. Math. 207 (2002), 257--286.

\bibitem{Brenti-2}
F. Brenti, Parabolic Kazhdan-Lusztig  $R$-polynomials for  tight quotients of the symmetric
group, J. Algebra 347 (2011), 247--261.




\bibitem{Deodhar}
V.V. Deodhar, On some geometric aspects of Bruhat orderings II,
The parabolic analogue of Kazhdan-Lusztig polynomials,
J. Algebra 111 (1987), 483--506.

\bibitem{Deodhar-2}
V.V. Deodhar, Duality in parabolic setup for questions in  Kazhdan-Lusztig theory,
J. Algebra 142 (1991), 201--209.


\bibitem{Humphreys}
J.E. Humphreys, Reflection groups and Coxeter groups,  Cambridge
Studies in Advanced Mathematics, No. 29, Cambridge Univ. Press,
Cambridge, 1990.


\bibitem{KL}
D. Kazhdan and G. Lusztig, Representations of Coxeter groups and
Hecke algebras, Invent. Math. 53 (1979), 165--184.

\bibitem{Stembridge} J. Stembridge, Tight quotients and double quotients in the Bruhat order,
Electron. J. Combin. 11 (2005), R14.

\end{thebibliography}
\end{document}